\documentclass[11pt,a4paper]{article}

\usepackage{mathtools,amssymb,amsthm,amsmath}
\usepackage[left=1.1in,right=1.1in]{geometry}
\usepackage{cite}

\usepackage{caption}
\usepackage{float}
\usepackage{subfloat}
\usepackage{url}
\usepackage{xcolor}
\usepackage{algorithmic}
\usepackage[linesnumbered,algoruled,boxed,lined]{algorithm2e}

\def\RR{\mathbb R}
\def\EE{\mathcal F}

\def\e{\varepsilon}

\def\argmin{{\rm arg}\!\min}
\def\be{\begin{equation}}
\def\ee{\end{equation}}
\def\bea{\begin{eqnarray}}
\def\eea{\end{eqnarray}}
\renewcommand{\phi}{\varphi}

\newcommand{\mc}[1]{\mathcal{#1}}

\newtheorem{remark}{Remark}[section]
\newtheorem{lemma}{Lemma}[section]
\newtheorem{theorem}{Theorem}[section]

\title{Optimization by linear kinetic equations\\ and mean-field Langevin dynamics}
\date{}

\author{Lorenzo Pareschi\thanks{Maxwell Institute for Mathematical Sciences and
Department of Mathematics, School of Mathematical and Computer Sciences (MACS), HWU Edinburgh, UK.} $^{,}$\thanks{Department of Mathematics and Computer Science, University of Ferrara, Italy.}}

\begin{document}

\maketitle

\begin{abstract}
Probably one of the most striking examples of the close connections between global optimization processes and statistical physics is the simulated annealing method, inspired by the famous Monte Carlo algorithm devised by Metropolis et al. in the middle of the last century. 
In this paper we show how the tools of linear kinetic theory allow to describe this gradient-free algorithm from the perspective of statistical physics and how convergence to the global minimum can be related to classical entropy inequalities. 
This analysis highlight the strong link between linear Boltzmann equations and stochastic optimization methods governed by Markov processes. Thanks to this formalism we can establish the connections between the simulated annealing process and the corresponding mean-field Langevin dynamics characterized by a stochastic gradient descent approach.
Generalizations to other selection strategies in simulated annealing that avoid the acceptance-rejection dynamic are also provided. 
\end{abstract}


\tableofcontents

\section{Introduction}
Since its introduction n the late 1980s^^>\cite{kirkpatrick1983optimization}, simulated annealing has become a popular optimization algorithm, and its applications have expanded to many fields, including artificial intelligence, machine learning, and operations research^^>\cite{3Dface10,HS05}. The algorithm was developed as an extension of the Metropolis-Hastings algorithm, a Monte Carlo method used for simulating complex systems in physics^^>\cite{metro53,Hast70}, and adapts the concept of annealing to optimization problems by viewing the process of slowly cooling a material as a search for the lowest-energy state of the system. 

Simulated annealing is similar to other metaheuristic algorithms such as genetic algorithms, ant colony optimization, particle swarm optimization and consensus based optimization in that it is based on the idea of exploring a large search space to find a global optimum without using gradient-based information^^>\cite{acta19}. However, simulated annealing is distinct in that it uses a single solution-based probabilistic approach to accept worse positions in the hope of finding a better one, whereas other metaheuristic algorithms often use a population-based approach and other stochastic strategies^^>\cite{Blum:2003:MCO:937503.937505, Gendreau:2010:HM:1941310, Back:1997:HEC:548530, pinnau2017consensus}.



By applying the above concepts, metaheuristic algorithms have been able to make significant advances in the search for valuable solutions to challenging optimization problems out of reach of traditional (gradient-based) methods. However, proving the rigorous convergence of metaheuristic optimization algorithms to the global minimum for non convex functionals, or to some reasonable approximation of it, remains a challenge. Indeed, metaheuristics 
involves the creative use of available resources to find efficient solutions without necessarily relying on a rigorous mathematical foundation that provides an analytical setting.


On the other hand, metaheuristics share similarities with statistical physics since they both deal with the complexities of large systems^^>\cite{Kolo:2010,Bellomo17}. The principles of statistical physics are versatile and powerful, providing insight into the behavior of large systems in a wide range of fields, from materials science to biophysics. By drawing upon the principles of statistical physics, it may be possible to provide a solid mathematical foundation to these class of optimization methods and develop more effective and efficient optimization algorithms that can handle increasingly complex problems and larger search spaces.

Mean field equations and kinetic equations are among the concepts in statistical physics that are most relevant to optimization. Mean field equations describe how each particle in a system interacts with a theoretical "average" field created by all the other particles in the system. This provides insight into the behavior of large systems, making it possible to predict macroscopic properties such as temperature or pressure. Kinetic equations, on the other hand, describe the evolution of a particle system, considering the interactions between particles as instantaneous, microscopic collisions^^>\cite{Villani02}.
These ideas have led to a new view of methaeuristic optimization by considering the corresponding continuous dynamics described by appropriate kinetic equations of Boltzmann type^^>\cite{Borghi1} and mean-field type^^>\cite{carrillo2018analytical, carrillo2019consensus,pinnau2017consensus,Kalise,Ha,FKR} even in constrained contexts^^>\cite{Borghi2,fhps1,fhps2,fhps3}, multiobjective situations^^>\cite{multiobj}, or in generalizations to sampling^^>\cite{Franca}. See^^>\cite{Claudia} for a recent survey.

Recently, leveraging these concepts has addressed the problem of providing rigorous theoretical support for some of the main optimization algorithms based on meta-heuristics, including particle swarm optimization^^>\cite{Grassi21, Hui23} and genetic algorithms^^>\cite{genetic,AFT}. In this paper we will focus our attention on one of the most notable examples of metaheuristics, namely the simulated annealing algorithm. This algorithm was inspired by the Monte Carlo algorithm developed by Metropolis et al. in the middle of last century. We show how classical tools of kinetic theory can be used to describe the Markov process which characterizes the method and show how its convergence to the global minimum is related to classical functional inequalities based on the so-called entropy method^^>\cite{Clem,Bisi15,Canizo18,DV01}. 

Although, from qualitative point of view, some of what is presented here has similarities to known results\cite{Aarts:1989:SAB:61990,holley1988simulated, Belisle92,Hajek}, nonetheless, this new perspective to the problem has several features that make it interesting. In the first place, it highlights the strong link between linear Boltzmann equations and stochastic optimization methods governed by Markov processes. Secondly, thanks to the novel formalism it is possible to establish the connections to the corresponding Langevin dynamics characterized by a mean-field description of a stochastic gradient descent method^^>\cite{Pavlio23,Chizat22,Mon18,GH86}. The latter result, which provides a way to establish a theoretical connection between gradient-based and gradient-free techniques,
 is obtained thanks to a quasi-invariant scaling limit inspired by the so-called grazing collisions asymptotics of the Boltzmann equation^^>\cite{partos13,DV92,grazing22}.


The rest of the manuscript is organized as described below. The next section introduces the classical simulated annealing algorithm and the corresponding Langevin dynamics. Section 3 is the devoted to the derivation of the kinetic description and the study of its long time behavior using classical functional inequalities. In particular, we show how the Langevin dynamics can be recovered in a suitable scaling limit and introduce a novel formulation of simulated annealing that avoid the acceptance-rejection dynamic. A few simple numerical experiments to validate the theoretical findings are then given in Section 4. Some concluding remarks are reported at the end of the manuscript.

\section{Simulated annealing}
In the sequel, to fix notations, we consider the following optimization problem
\begin{equation}\label{typrob}
x^\ast \in \argmin\limits_{x\in \RR}\EE(x)\,,
\end{equation}
where $\EE(x):\mathbb R^{d} \to \mathbb R$ is a given high-dimensional, non convex, cost function, which we wish to minimize. We will assume that the minimizing argument $x^*\in\RR^d$ of \eqref{typrob} exists and is unique.
This kind of problems are ubiquitous in many applications, most recently and notably in {machine learning} (e.g., training of neural networks), {signal/image processing} and optimal control of multi-agent systems^^>\cite{carrillo2019consensus,fhps2,Borghi2,bishop06}.

\subsection{The metaheuristic algorithm}
\label{sec:sa}
Simulated annealing is a metaheuristic optimization algorithm that was developed by Kirkpatrick, Gelatt, and Vecchi in 1983^^>\cite{kirkpatrick1983optimization}. The algorithm was inspired by the Monte Carlo algorithm of Metropolis et al.^^>\cite{metro53}, dating from 30 years earlier, simulating the thermal motion of particles in thermal contact with a heat bath at a certain temperature. As it is well-known the method permits to recover the Boltzmann-Gibbs distribution of particles. 

More precisely, the idea behind simulated annealing is to start with an initial solution and gradually modify it by accepting worse solutions with a certain probability, allowing the algorithm to escape local optima and explore different regions of the solution space. The probability of accepting a worse solution is controlled by a temperature parameter that is gradually reduced over time, mimicking the cooling process of annealing in metallurgy from which it inherited the name.

\begin{figure}
\begin{center}
\includegraphics[scale=0.4]{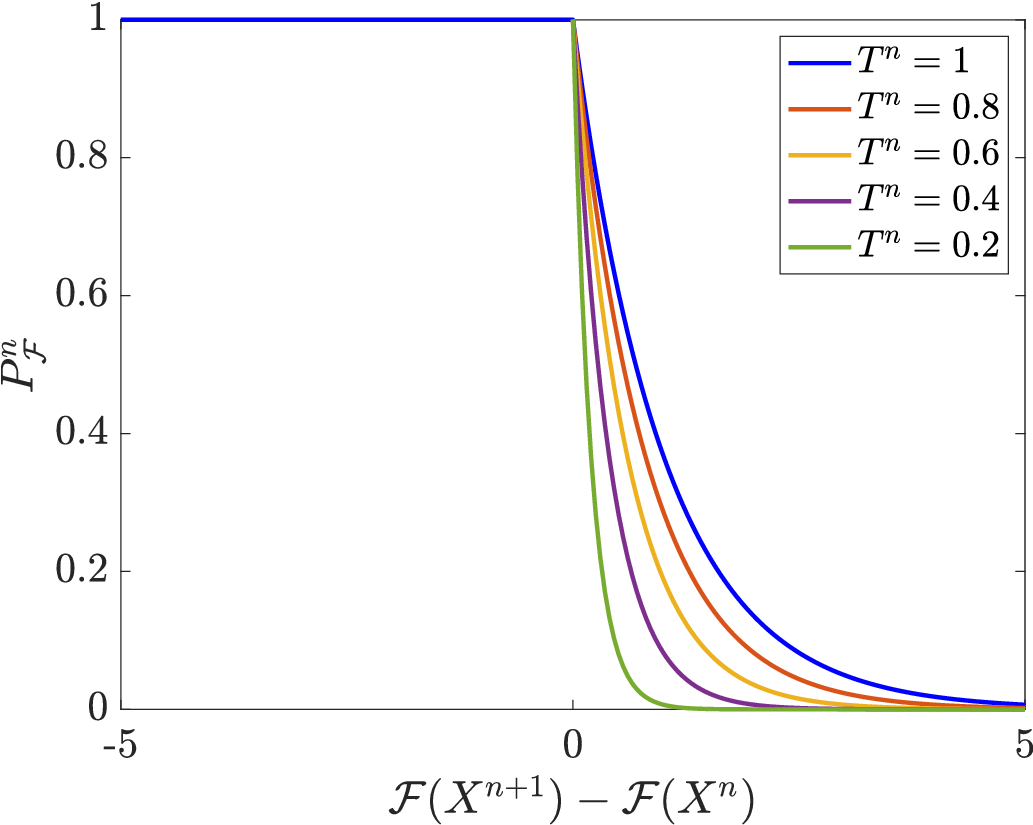}
\end{center}
\caption{The probability to accept a trial point in simulated annealing}
\label{fg:figp}
\end{figure}

One step of the standard simulated annealing algorithm, starting from a random trial point $X^0$ and an initial control parameter $T^0$, referred to as \emph{temperature}, can be summarized for $n\geq 0$ in the following steps.
\begin{enumerate}
\item Choose the distance of the new trial point from the current point by a probability distribution with zero mean and a time dependent variance 
\be
\tilde X^{n+1} = X^{n} + \sigma^n\xi
\label{s1}
\ee
where $\xi \sim p(\xi)$, where, for example, $p(\xi)$ is a uniform or a normal probability density with mean zero and unitary variance, and $\sigma^n > 0$ depends on the temperature $T^n$. A classical choice is $\sigma^n$ proportional to $\sqrt{T^n}$ so that the temperature is related to the variance of the random noise.
\item The algorithm determines whether the new point $\tilde X^{n+1}$ is better or worse than the current point $X^n$. If the new point is better than the current point, it becomes the next point. If the new point is worse than the current point, the algorithm can still accept it in accordance with a temperature dependent Boltzmann-Gibbs probability. In summary 
\be
X^{n+1} = \begin{cases}
\tilde X^{n+1} & {\rm if}\,\EE(\tilde X^{n+1})-\EE(X^n) < 0, \\
\tilde X^{n+1} \quad\hbox{with probability}\,\, e^{-\frac{\EE(\tilde X^{n+1})-\EE(X^n)}{T^n}} & \\[-.3cm]
& {\rm if}\,\EE(\tilde X^{n+1})-\EE(X^n) \geq 0,\\[-.2cm]
 X^n \qquad\hbox{otherwise.} & 
\end{cases}
\label{s2}
\ee
\item The algorithm systematically lowers the temperature, accordingly to a law of the type
\be
T^{n+1}=\lambda^{n} T_0,\qquad \lambda^n \in (0,1),
\label{s3}
\ee
where $T_0>0$ is a given initial temperature. A classical choice is $\lambda^n = 1/\ln(n+2)$ for which the algorithm has been shown to converge in the discrete case where the search space has a finite sets of states^^>\cite{Hajek,Geman84}. 
\end{enumerate}
We can rewrite points 1 and 2 of the above process in compact form as follows
\be
\begin{split}
X^{n+1} &= \Psi\left(\eta <  P_\EE^n\right)(X^{n} + \sigma^n\xi) +\left(1-\Psi\left(\eta <  P_\EE^n\right)\right)X^n\\
&= X^n + \sigma^n\Psi\left(\eta <  P_\EE^n\right)\xi
\label{sacom}
\end{split}
\ee
with $\eta \sim U(0,1)$, $\Psi(\cdot)$ the indicator function and we introduced the probability to accept a trial point $P_\EE^n$  (see Figure \ref{fg:figp}) 
\be
P_\EE^n =
\begin{cases}
1, & \EE(X^{n+1}) < \EE(X^{n}),\\
e^{-\frac{\EE(X^{n+1})-\EE(X^{n})}{T(t)}}, & \EE(X^{n+1}) \geq \EE(X^{n}).
\end{cases}
\label{pro}
\ee
The algorithm just described is usually referred to as classical simulated annealing. 
Since its introduction, several variations and extensions of the algorithm have been proposed, including parallel simulated annealing, adaptive simulated annealing, and hybrid simulated annealing with other optimization methods^^>\cite{Aarts:1989:SAB:61990,holley1988simulated,TSALLIS1996395}.
Various convergence results have also been published in the literature, see for example^^>\cite{Belisle92, Locatelli00}.

\subsection{Mean-field Langevin dynamics}
A continuous version of simulated annealing has been proposed in^^>\cite{GH86} described by a stochastic differential process in the form
\be
dX_t = -\nabla_x \EE(X_t)dt+\sqrt{2T}dB_t,
\label{eq:csa}
\ee 
where $B_t$ is a standard Brownian motion. The convergence analysis for this dynamic, sometimes referred to as \emph{Langevin equation} or \emph{overdamped kinetic simulated annealing}, is much more rich and has been considered in several recent papers, we refer to^^>\cite{Pavlio23, Mon18} and the references therein. Note that the stochastic differential system \eqref{eq:csa} can likewise be understood as the limit for small learning rates of a stochastic gradient descent (SGD) method.

The reason why process \eqref{eq:csa} is refereed to as continuous simulated annealing is because its mean field description^^>\cite{Chizat22, HRSS21}
\be
\frac{\partial f}{\partial t}(x,t) = \nabla_x\cdot \left(\nabla_x \EE(x) f(x,t)\right) + {T}\Delta_{xx} f(x,t),
\label{eq:mfsa}
\ee
where $f(x,t)$ is the probability density to have a trial point in position $x\in\RR^d$ at time $t>0$,
under suitable conditions on $\EE(x)$, admits as stationary state the Boltzmann-Gibbs distribution
\be
f^\infty_\EE(x)=C e^{\frac{-\EE(x)}{T}},
\label{eq:Gibbs}
\ee
with $C$ a normalization constant. For small values of the temperature $T$, the stationary state $f^\infty_\EE(x)$ concentrates on global minima of $\EE(x)$ and we expect to find $X_t$ in \eqref{eq:csa} near a global minimum. Unfortunately, the time to reach equilibrium increases exponentially with $1/T$, so that solutions to \eqref{eq:csa} or \eqref{eq:mfsa} with small $T$ will be extremely slow to find relevant minima^^>\cite{GH86}. 

Similar to simulated annealing, this suggest slowly decreasing the temperature with time during the dynamic, so that the solution can approach the equilibrium at a faster rate and concentrate on minima asymptotically. More precisely, one can show that, as $T(t)= T_0/\log(2+t)$ for $T_0$ sufficiently large, the above probability
measure converges weakly to the set of global minima based on the Laplace principle^^>\cite{Hwang}, which states that, for any absolutely continuous probability distribution $g(x)$ on $\RR^d$, we have
\be
\lim_{T\to 0} -T\log\left(\int_{\RR^d} g(x)e^{\frac{-\EE(x)}{T}}\,dx\right) = \inf_{x\in{\rm supp}(g)} \EE(x).
\label{eq:Laplace}
\ee 
It should be pointed out that for a time-dependent control temperature, the normalization constant $C$ and consequently the solution \eqref{eq:Gibbs} depend on time.

However, the diffusion process \eqref{eq:csa} was simply inspired by the original simulated annealing algorithm, but was not derived directly from the discrete algorithm \eqref{s1}-\eqref{s3} since, in the words of the authors of^^>\cite{GH86}, ``the extension of the Metropolis algorithm to continuous variables involves some awkward computational problems". Probably, the main limitation of \eqref{eq:csa} is that it requires the gradient evaluation, in contrast with the gradient-free meta-heuristic nature of the Metropolis algorithm. In the following we will adopt a different approach to derive a gradient-free continuous extension of the simulated annealing algorithm \eqref{s1}-\eqref{s3} and we will try to shed some light on the relationships between the original metaheuristic algorithm and the stochastic differential process \eqref{eq:csa}.



\section{Linear kinetic equations for global optimization}
A continuous model for the dynamic described by \eqref{sacom}-\eqref{pro} can be derived using the tools of linear kinetic theory and observing that \eqref{pro} represents the probability to accept the post-interaction position. We remark that the simulated annealing dynamics is often studied in the framework of Markov chain Monte Carlo (MCMC) processes^^>\cite{Hast70}. Here we follow a different, although closely related, approach based on linear Boltzmann-type equations^^>\cite{Kolo:2010, partos13, Bisi15, Canizo18}. The new formalism permits to link convergence to the global minimum to classical functional inequalities based on the entropy method^^>\cite{Clem,Bisi15,Canizo18,DV01} and to establish the connections with the corresponding mean-field Langevin dynamics thanks to a suitable scaling limit inspired by the so-called grazing collisions asymptotics of the Boltzmann equation^^>\cite{partos13,DV92,grazing22}. 

\subsection{Kinetic description of simulated annealing}
After introducing the probability density $f(x,t)$, we can write the evolution equation 
\be
\begin{split}
\frac{\partial f(x,t)}{\partial t} &= \left\langle B_\EE(x'\to x)f(x',t) - B_\EE(x\to x') f(x,t)\right\rangle:={\mathcal L}_\EE(f(x,t))
\end{split}
\label{eq:ksa}
\ee
where we used the notation 
\[
\langle\,g\,\rangle=\mathbb{E}_\xi[g] = \int_{\RR^d} g(\xi)p(\xi)d\xi,
\]
to denote the expectation of a function $g(\xi)$ with respect to the {selection probability} $p(\xi)$, $\xi\in \RR^d$. Furthermore, $x'$ is the new trial-point position given by
\be
x'=x+\sigma(t)\xi,
\ee
and 
\be
B_\EE(x\to x')=\min\left\{1,\frac{f_\EE^\infty(x')}{f_\EE^\infty(x)}\right\}.
\label{eq:kernel}
\ee 
The first term on the r.h.s of equation \eqref{eq:ksa} characterizes the expected gain of particles in position $x$ starting from position $x'$ with transition probability $B_\EE(x'\to x)$, whereas the second term the expected loss of particles in position $x$ because of a move to position $x'$ with transition probability $B_\EE(x\to x')$. In the sequel, we will assume $p(\xi)$ symmetric, i.e. $p(-\xi)=p(\xi)$, with mean $0$ and identity covariance matrix $\Sigma=I_d$.

We also assume a given cooling law $T(t)$ such that $dT/dt \leq 0$ for the decay of the temperature. The resulting linear kinetic model \eqref{eq:ksa} belongs to the general class of linear kinetic models, see also^^>\cite{partos13,Furioli20,Canizo18,Bisi15,Russ}.  

Let us remark that since $\left\langle B_\EE(x'\to x) \right\rangle \leq 1$
the loss term in \eqref{eq:ksa} is bounded and then the operator $\mathcal J_\EE(f) = {\mathcal L}_\EE(f)+f$ is a positive monotone operator in the sense that
\be
\mathcal J_\EE(f) \geq  \mathcal J_\EE(g)\geq 0\qquad \hbox{if}\qquad f\ge g \geq 0.
\label{eq:J}
\ee
Additionally by mass conservation we have $\int_{\RR^d} \mathcal J_\EE(f)\,dx =1$ so that $\mathcal J_\EE(f)$ is a probability density function.

Using Young's inequality we get
\be
\begin{split}
\|\mathcal J_\EE(f)\|_{L^p(\RR^d)} \leq \|p\ast f\|_{L^p(\RR^d)}+(1-\lambda_1)\|f\|_{L^p(\RR^d)}
& \leq (2-\lambda_1)\|f\|_{L^p(\RR^d)}
\end{split}
\ee
where $\lambda_1=e^{\frac{-\Delta}{T}}$, $\Delta = \sup_{x,y\in\RR^d}\left\{\EE(x)-\EE(y)\right\}$.
From the above bounds, by standard arguments, one obtains existence and uniqueness of a nonnegative solution for \eqref{eq:ksa} with nonnegative initial data.

The following simple identity will be useful in the sequel
\begin{lemma}
For any symmetric probability density $p(\xi)$ and any integrable function $g(x,x')$ we have
\be
\left\langle \int_{\RR^{d}} g(x',x)\,dx\right\rangle=\left\langle \int_{\RR^{d}}g(x,x')\,dx\right\rangle.
\label{eq:ident}
\ee
\label{le:change}
\end{lemma} 
\proof
The identity is obtained by the change of integration variables $(\xi,x)\to (-\xi,x')$ and then switching notations back from $x'\to x$. In fact, by the simmetry of $p(\xi)$ we get
\[
\left\langle \int_{\RR^{d}} g(x+\sigma\xi,x)\,dx\right\rangle=\left\langle \int_{\RR^{d}} g(x',x'+\sigma\xi)\,dx'\right\rangle=\left\langle \int_{\RR^{d}} g(x,x+\sigma\xi)\,dx\right\rangle.
\]
\endproof
Thanks to the above identity, system \eqref{eq:ksa} can be conveniently written in weak form as
\be
\frac{\partial}{\partial t}\int_{\RR^d} f(x,t)\phi(x)\,dx = \left\langle \int_{\RR^{d}}B_\EE(x\to x')(\phi(x')-\phi(x))f(x,t)\,dx\right\rangle,
\label{weaksa}
\ee
for any given sufficiently regular function $\phi=\phi(x)$. 

Clearly when $\phi(x)=1$ we get conservation of the total number of trial points, namely starting from a probability density we have a probability density at all later times. Using a kinetic theory terminology, this is the unique collision invariant of the system. In fact, one cal easily check that the expected value is not conserved since it is driven by a combination of a descent dynamic with a low probability ascent one. Analogous arguments applies also to the variance.

Furthermore, thanks again to \eqref{eq:ident} we also have the following weak form
\be
\begin{split}
\frac{\partial}{\partial t}\int_{\RR^d} f(x,t)\phi(x)\,dx &=\\
&\hskip -1cm \frac12\left\langle \int_{\RR^{d}}(\phi(x')-\phi(x))(B_\EE(x\to x')f(x,t)-B_\EE(x'\to x)f(x',t))\,dx\right\rangle.
\end{split}
\ee
It is immediate to verify that  
\begin{lemma}
The Boltzmann-Gibbs distribution \eqref{eq:Gibbs} satisfies ${\mathcal L}_{\EE}(f^\infty_\EE(x))=0$, $\forall\,x\in\RR^d$.
\end{lemma}
\proof
It is enough to observe that we have the detailed balance condition
\be
B_\EE(x'\to x)f^\infty_\EE(x')-B_\EE(x\to x')f^\infty_\EE(x) = 0.
\label{eq:db}
\ee 
In fact, by direct substitution for $\EE(x')<\EE(x)$ we obtain
\[
\frac{f^\infty_\EE(x)}{f^\infty_\EE(x')}f^\infty_\EE(x') -f^\infty_\EE(x) = 0,
\]
whereas for $\EE(x')>\EE(x)$ we get
\[
f^\infty_\EE(x')-\frac{f^\infty_\EE(x')}{f^\infty_\EE(x)}f^\infty_\EE(x) = 0.
\]
\endproof
Note that here, unlike in classical kinetic theory, it is the transition probabilities and not the microscopic interactions that play an essential role in determining the equilibrium state.
\begin{remark} 
\label{rem1}
The equation \eqref{eq:ksa}, or its weak form \eqref{weaksa}, admit other equivalent formulations which emphasize connections with stochastic processes and classical kinetic theory. 
\begin{itemize}
\item
In fact, equation \eqref{eq:ksa} represents the Kolmogorov forward equation associated with a Markov process on $\RR^d$ described by the simulated annealing algorithm with invariant measure $f^\infty_\EE(x)$. If we define the nonnegative kernel of the operator ${\mathcal L}_\EE(f(x,t))$ as
\be
k_\EE(x',x)=\frac1{\sigma(t)^d}p\left(\frac{x'-x}{\sigma(t)}\right) B_\EE(x'\to x)
\label{eq:kF}
\ee 
we can write 
\be
{\mathcal L}_\EE(f(x,t))= \int_{\RR^d} k_\EE(x',x) f(x')\,dx'-f(x,t)\int_{\RR^d} k_\EE(x,x')\,dx'.
\label{eq:carl}
\ee
\item
It is possible to write the weak form \eqref{weaksa}  emphasizing the similarities with classical linear kinetic equations of Boltzmann type as
\be
\frac{\partial}{\partial t}\int_{\RR^d} f(x,t)\phi(x)\,dx = \int_{\RR^d}\int_{\RR^{d}}\beta_\EE(x\to x')(\phi(x')-\phi(x))f(x,t)f^\infty_\EE(x')\,dx\,dx',
\label{weaksa2}
\ee
where $\beta_\EE(x\to x') \geq 0$ is now a symmetric collision kernel
\be
\beta_\EE(x\to x') = \frac1{\sigma(t)^d}p\left(\frac{x'-x}{\sigma(t)}\right)\min\left\{\frac1{f^\infty_\EE(x')},\frac1{f^\infty_\EE(x)}\right\}.
\ee
In this latter formulation, the stationary solution is $f_\EE^\infty(x)$ regardless of the nonnegative collision kernel $\beta_\EE(x\to x')$. For example, the simplified choice $\beta_\EE(x\to x')=1$ 
yields
\be
{\mathcal L}_\EE(f(x,t))= \left\langle f(x',t)\right\rangle f^\infty_\EE(x)-\left\langle f_\EE^\infty(x')\right\rangle f(x,t).
\label{eq:Maxwell}
\ee
From an algorithmic perspective, the dynamic described by \eqref{eq:Maxwell} corresponds to a random walk where we accept each new move $x'$ with probability given by the Boltzmann-Gibbs measure $f^\infty_\EE(x')$. 


\end{itemize}
\end{remark}

\subsection{Convergence to equilibrium}
An important result to study the convergence to equilibrium of linear Boltzmann equations is to observe that taking 
\be
\phi(x)=\log\left(\frac{f(x,t)}{f^\infty_\EE(x,t)}\right),
\ee
one can prove the decay of the so-called \emph{relative Shannon-Boltzmann entropy} (also known as \emph{Kullback--Leibler divergence}) along solutions of equation \eqref{weaksa}. We refer to^^>\cite{Canizo18, Bisi15, Perthame05} for this type of results in the case of other linear kinetic equations. Let us notice that there are two main differences here. Firstly, the kernel \eqref{eq:kernel} is not symmetric due to the optimization process, which means that moments are not preserved except for the total number of trial points. Secondly, the kernel and has a fundamental rule in characterizing the long time behavior of the equation, which, as we will discuss later, may depend on time due to the control induced by the annealing temperature.

Let us first consider the case with constant temperature, so that $C$ in \eqref{eq:Gibbs} is a constant. We have 
\[
\begin{split}
\frac{\partial}{\partial t}\int_{\RR^d} f(x,t)\log\left(\frac{f(x,t)}{f^\infty_\EE(x,t)}\right)\,dx 
&= \frac12\left\langle \int_{\RR^{d}}\left(\log\left(\frac{f(x',t)}{f^\infty_\EE(x')}\right)-\log\left(\frac{f(x,t)}{f^\infty_\EE(x)}\right)\right)\right.\\
&\qquad (B_\EE(x\to x')f(x,t)-B_\EE(x'\to x)f(x',t))\,dx\Big\rangle
\end{split}
\]
Now, when $\EE(x') < \EE(x)$ we get
\[
B_\EE(x\to x')f(x,t)-B_\EE(x'\to x)f(x',t) = f(x,t)-\frac{f^\infty_\EE(x)}{f^\infty_\EE(x')}f(x',t)
\]
so that the integrand becomes
\be
f^\infty_\EE(x)\left(\log\left(\frac{f(x',t)}{f^\infty_\EE(x')}\right)-\log\left(\frac{f(x,t)}{f^\infty_\EE(x)}\right)\right)\left(\frac{f(x,t)}{f^\infty_\EE(x)}-\frac{f(x',t)}{f^\infty_\EE(x')}\right)\leq 0,
\ee
since the function 
\be
h(x,y)=(x-y)(\log(x)-\log(y)) \geq 0,
\label{eq:h}
\ee
for all $x,y\in\RR^+$.

Similarly, when $\EE(x) < \EE(x')$ we can write
\[
B_\EE(x\to x')f(x,t)-B_\EE(x'\to x)f(x',t) = \frac{f^\infty_\EE(x')}{f^\infty_\EE(x)}f(x,t)-f(x',t),
\]
and the integrand now satisfies
\[
f^\infty_\EE(x')\left(\log\left(\frac{f(x',t)}{f^\infty_\EE(x')}\right)-\log\left(\frac{f(x,t)}{f^\infty_\EE(x)}\right)\right)\left(\frac{f(x,t)}{f^\infty_\EE(x)}-\frac{f(x',t)}{f^\infty_\EE(x')}\right)\leq 0.
\]
As a consequence we have the following result
\begin{theorem}
The functional
\be
H(f|f^\infty_\EE)=\int_{\RR^d} {f(x,t)}\log\left(\frac{f(x,t)}{f^\infty_\EE(x)}\right)\,dx,
\label{eq:ent}
\ee
satisfies
\be
\frac{dH(f|f^\infty_\EE)}{dt} = -I_\EE[f]\leq 0,
\label{eq:ineq}
\ee
where
\be
I_\EE[f] = \frac12\left\langle \int_{\RR^d}B_\EE(x\to x')f^\infty_\EE(x)\,h\left(\frac{f(x',t)}{f^\infty_\EE(x')},\frac{f(x,t)}{f^\infty_\EE(x)}\right)\,dx\right\rangle,
\label{eq:ediss}
\ee
and the function $h$ is given by \eqref{eq:h}.
\end{theorem}
\proof
By direct computation we get
\[
\begin{split}
\frac{dH(f|f^\infty_\EE)}{dt} &= \int_{\RR^d} \frac{\partial f(x,t)}{\partial t}\left(1+\log\left(\frac{f(x,t)}{f^\infty_\EE(x)}\right)\right)\,dx\\
& = \int_{\RR^d} \frac{\partial f(x,t)}{\partial t}\log\left(\frac{f(x,t)}{f^\infty_\EE(x)}\right)\,dx =- I_\EE[f]\leq 0.
\end{split}
\]
\endproof
\noindent
\begin{remark}
The above arguments can be extended to any convex function $\Phi(x)$ so that 
\be
H_\Phi(f|f^\infty_\EE)=\int_{\RR^d} {f^\infty_\EE(x)}\Phi\left(\frac{f(x,t)}{f^\infty_\EE(x)}\right)\,dx,
\label{eq:hphi}
\ee
is a Lyapunov functional satisfying \eqref{eq:ineq}-\eqref{eq:ediss} for $h(x,y)= (x-y)(\Phi'(x)-\Phi'(y))\geq 0$. In the particular case $\Phi(x)=x\log(x)-x+1$ the expression \eqref{eq:hphi} reduces to the Boltzmann-Shannon entropy \eqref{eq:ent}, whereas for $\Phi(x)=(x-1)^2/2$ the corresponding dissipation \eqref{eq:ediss} is represented by the Dirichlet form
\be
D_\EE[f]=\frac12 \left\langle \int_{\RR^d}B_\EE(x\to x')f^\infty_\EE(x)\,\left(\frac{f(x',t)}{f^\infty_\EE(x')}-\frac{f(x,t)}{f^\infty_\EE(x)}\right)^2\,dx\right\rangle.
\label{eq:Dirichlet}
\ee
Note also that, using \eqref{eq:kF}, we can write  \eqref{eq:ediss} as
\be
I_\EE[f] = \frac12\int_{\RR^d} \int_{\RR^d}k_\EE(x,x') f^\infty_\EE(x)\,h\left(\frac{f(x',t)}{f^\infty_\EE(x')},\frac{f(x,t)}{f^\infty_\EE(x)}\right)\,dx\,dx'.
\label{eq:ediss2}
\ee

See for example^^>\cite{Furioli20, Bobko06, holley1988simulated,Clem} and the references therein for other convex functionals that can be used to study the large-time behavior of the solution to equation \eqref{eq:ksa}.
\end{remark}
It is well-known that if we can find an estimate of the type^^>\cite{DV01}
\be
I_\EE[f] \geq \lambda H(f|f^\infty_\EE)
\label{eq:ccg}
\ee
with $\lambda > 0$ we get for an initial data $f(x,0)=f_0(x)$
\[
H(f|f^\infty_\EE) \leq H(f_0|f^\infty_\EE)e^{-\lambda t}
\]
and, due to the Csisz\'ar–Kullback inequality 
\[
H(f|f^\infty_\EE) \geq \frac{\|f-f_\EE^\infty\|_{L_1(\RR^d)}}{2\|f\|_{L_1(\RR^d)}},
\]
the convergence of $H(f|f^\infty_\EE)$ towards $0$ implies the convergence in $L_1(\RR^d)$ of $f(x,t)$ towards $f^\infty_\EE(x)$. Therefore, the Shannon-Boltzmann entropy decays until it reaches its unique minimum $f^\infty_\EE(x)$. Analogous estimates for the Dirichlet functional \eqref{eq:Dirichlet} characterize the spectral gap of the equation and will lead to decays in the $L^2(\RR^d)$ norm with weight ${f^\infty_\EE}^{-1}$. Estimates of type \eqref{eq:ccg} for the Dirichlet functional are known in the literature for simulating annealing type dynamics on a finite state space^^>\cite{holley1988simulated} and to the best of the author's knowledge are not known for the continuous formulation here presented.

In the general case where the temperature is time dependent we must take into account the normalization constant and consider that
\[\phi(x)=\log\left(\frac{f(x,t)}{f^\infty_\EE(x,t)}\right)=\log(f(x,t))+\frac{\EE(x)}{T(t)}+\log\left(C(t)\right)
\] 
to get
\be
\begin{split}
&\frac{d}{dt}\int_{\RR^d} f(x,t)\left(\log(f(x,t))+\frac{\EE(x)}{T(t)}+\log\left(C(t)\right)\right)\,dx\\
&\quad=\int_{\RR^d}\frac{\partial f(x,t)}{\partial t} \left(\log(f(x,t))+\frac{\EE(x)}{T(t)}\right)\,dx\\
&\qquad  
-\frac{T'(t)}{T^2(t)}\int_{\RR^d} f(x,t)\EE(x)\,dx + \frac{C'(t)}{C(t)}.
\end{split}
\ee
By direct differentiation of the normalization constant we get
\[
\frac{C'(t)}{C(t)} = \frac{T'(t)}{T^2(t)} \int_{\RR^d} f^\infty_\EE(x,t)\EE(x)\,dx, 
\]
so that we need to control the additional term
\[
-\frac{T'(t)}{T^2(t)}\int_{\RR^d} \EE(x)\left(f(x,t)-f^\infty_\EE(x,t)\right)\,dx.
\]
This requires that $T'(t)=o(T^2(t))$ as $T(t)\to 0$. This shows that the temperature decay has to be slow enough. For example if $T(t) \approx 1/t$ we get $T'(t)/T(t)^2 \approx 1$ whereas for $T(t)\approx 1/\log(t)$ we get $T'(t)/T(t)^2 \approx 1/t$ and the quantity can be bounded.

Now, using again Cszizar-Kullback inequality we can write 
\[
\begin{split}
\frac{d\|f-f^\infty_\EE\|_{L_1(\RR^d)}}{dt} &\leq -\left(\lambda-\frac{c}{t} \|\EE\|_{L_\infty(\RR^d)}\right) \|f-f^\infty_\EE\|_{L_1(\RR^d)},
\end{split}
\]
where $c>0$ is a suitable constant.

Therefore, we can state the following result.

\begin{theorem} If the entropy dissipation rate satisfies \eqref{eq:ccg}, then for a time dependent temperature $T(t)=O(1/\log(t))$ we have the entropy dissipation 
\be
\frac{dH_\EE(f|f^\infty_\EE)}{dt} \leq -\tilde\lambda(t) H_\EE(f|f^\infty_\EE),
\label{eq:est}
\ee
where $\tilde\lambda(t)=\lambda-\frac{c}{t} \|\EE\|_{L_\infty(\RR^d)} > 0$ when $t > c \|\EE\|_{L_\infty(\RR^d)}/\lambda$.   
\end{theorem}
Since by the Laplace principle \eqref{eq:Laplace}, as $t\to\infty$, $T(t)\to 0$ we have that $f^\infty_\EE(x,t)$ concentrates on the global minimum $x^*$, then by \eqref{eq:est} also $f(x,t)$ concentrates on $x^*$ and thus convergence to the global minimum can be obtained by standard arguments^^>\cite{Hwang, FKR, Chizat22,GH86}.

\subsection{Mean-field limit and Langevin dynamics}
One may ask whether the mean field description \eqref{eq:mfsa} is related to the kinetic description \eqref{eq:ksa}. Before entering in to the discussion let us observe that the weak form \eqref{weaksa} of the kinetic equation can be rewritten as follows
\[
\begin{split}
\frac{\partial}{\partial t}\int_{\RR^d} f(x,t)\phi(x)\,dx &= \left\langle\int_{\RR^d} (\phi(x')-\phi(x))f(x,t)\,dx\right\rangle \\
&\quad - \left\langle \int_{\RR^d}\left(1-\frac{f^{\infty}_\EE(x')}{f^\infty_\EE(x)}\right)\Psi(\EE(x')\geq \EE(x))(\phi(x')-\phi(x))f(x,t)\,dx\right\rangle.
\end{split}
\]

By analogy with the grazing collision limit of the Boltzmann equation, a way to derive the mean-field approximation from the corresponding kinetic formulation is based on the well-established approach in^^>\cite{partos13}. 
This corresponds to consider the behavior of the system under the scaling
\be
t \to t/\varepsilon,\quad \sigma(t) \to \sqrt{\varepsilon} \sigma(t),
\label{eq:scal}
\ee
and to write for small values of $\varepsilon \ll 1$
\be
\phi(x')=\phi(x)+(x'-x)\cdot\nabla_x \phi(x)+\frac12 \sum_{i,j=1}^d (x'_i-x_i)(x'_j-x_j) \frac{\partial^2\phi(x)}{\partial x_i\partial x_j}+O(\varepsilon^{3/2}).
\label{eq:expand1}
\ee
Additionally, due to its dependence from the post-interaction position, we also expand the equilibrium density as
\be
\begin{split}
f^{\infty}_\EE(x') &= f^\infty_\EE(x) + (x'-x) \cdot \nabla_x f^\infty_\EE(x)+O(\varepsilon)\\
& = f^\infty_\EE(x) - (x'-x) \cdot \frac1{T(t)}(\nabla_x \EE(x)) f^\infty_\EE(x)+O(\varepsilon).
\end{split}
\label{eq:expand2}
\ee
Under the scaling \eqref{eq:scal}, using \eqref{eq:expand1}-\eqref{eq:expand2} and the fact that $\xi$ has zero mean, 
we get the approximated weak form up to $O(\sqrt{\varepsilon})$
\[
\begin{split}
\frac{\partial}{\partial t}\int_{\RR^d} f(x,t)\phi(x)\,dx &= \frac{\sigma(t)^2}{2}\sum_{i,j=1}^d\int_{\RR^d}
p(\xi) \xi_i \xi_j\,d\xi \int_{\RR^{d}} \frac{\partial^2\phi(x)}{\partial x_i\partial x_j} f(x,t)\,dx\\
&\quad - \frac{\sigma(t)^2}{T(t)}\left\langle \int_{\RR^d}\Psi(\EE(x')\geq \EE(x))\xi\cdot \nabla_x \EE(x)\xi\cdot\nabla_x\phi(x)f(x,t)\,dx\right\rangle.
\end{split}
\]
Now using the symmetry of $p(\xi)$ and Lemma \ref{le:change} we have
\[
\begin{split}
&\left\langle\int_{\RR^d}\Psi(\EE(x')\geq \EE(x))\xi\cdot \nabla_x \EE(x)\xi\cdot\nabla_x\phi(x)f(x,t)\,dx\right\rangle\\
&\quad =\left\langle\int_{\RR^d}\Psi(\EE(x)\geq \EE(x'))\xi\cdot \nabla_x \EE(x')\xi\cdot\nabla_x\phi(x')f(x',t)\,dx\right\rangle=I,
\end{split}
\]
which using again Taylor expansions on the integrands yields
\[
\quad I= \left\langle\int_{\RR^d}\Psi(\EE(x)\geq \EE(x'))\xi\cdot \nabla_x \EE(x)\xi\cdot\nabla_x\phi(x)f(x,t)\,dx\right\rangle+O(\sqrt{\varepsilon}).
\]
Thus, assuming that for each fixed $x$ the set of points $y\in\RR^d$ s.t. $\EE(y)=\EE(x)$ has zero measure we get
\[
\begin{split}
&\left\langle\int_{\RR^d}\Psi(\EE(x')\geq \EE(x))\xi\cdot \nabla_x \EE(x)\xi\cdot\nabla_x\phi(x)f(x,t)\,dx\right\rangle\\
&= \frac12 \left\langle\int_{\RR^d}\xi\cdot \nabla_x \EE(x)\xi\cdot\nabla_x\phi(x)f(x,t)\,dx\right\rangle+O(\sqrt{\varepsilon}).
\end{split}
\]
Since
\be
\int_{\RR^d}
p(\xi) \xi_i \xi_j\,d\xi = \delta_{ij},\forall\, i,j=1,\ldots,d
\label{eq:var}
\ee
where $\delta_{ij}$ is the Kronecker delta,
sending $\varepsilon\to 0$ provided $f$ has bounded third order derivatives and $\EE(x)$ bounded second order derivatives, we formally have up to $O(\sqrt{\e})$
\be
\begin{split}
\frac{\partial}{\partial t}\int_{\RR^d} f(x,t)\phi(x)\,dx &= \frac{\sigma(t)^2}{2}\sum_{i=1}^d \int_{\RR^{d}} \frac{\partial^2\phi(x)}{\partial x_i^2} f(x,t)\,dx\\
&\quad - \frac{\sigma(t)^2}{2T(t)}\left\langle\int_{\RR^d}\xi\cdot \nabla_x \EE(x)\xi\cdot\nabla_x\phi(x)f(x,t)\,dx\right\rangle.
\end{split}
\label{eq:exp1}
\ee
Taking $2T(t)=\sigma^2(t)$, using again \eqref{eq:var}, we can revert to the original variables to get the mean-field Langevin dynamic
\be
\frac{\partial f(x,t)}{\partial t}=\nabla_x\cdot\left(\nabla_x\EE(x)f(x,t)\right)+T(t)\Delta_{xx} f(x,t).
\label{eq:VlasovSA}
\ee

\subsection{Maxwellian simulated annealing}
In this section we present an interesting alternative formulation of simulated annealing type processes which shares the common feature of originating the same mean field Langevin dynamic \eqref{eq:VlasovSA}. This formulation, however, at the kinetic level is more difficult to analyze and do not allow to explicitly compute the steady state solution. Therefore, we limit our considerations to the behavior in the mean field scaling \eqref{eq:scal} leaving a further discussion to the numerical section.

We can formulate the simulated annealing process avoiding the acceptance rejection dynamic by following classical analogies with kinetic equation with constant interaction kernels, usually refereed to as Maxwellian kernels^^>\cite{partos13}. To this aim starting from the trial point 
\be
\tilde X^{n+1} = X^n + \sigma^n \xi,
\ee
we define
\be
X^{n+1} = \begin{cases}
\tilde X^{n+1} & {\rm if}\,\EE(\tilde X^{n+1})-\EE(X^n) < 0 \\
X^n + e^{-\frac{\EE(\tilde X^{n+1})-\EE(X^n)}{T^n}}(\tilde X^{n+1} -X^n) & 
 {\rm if}\,\EE(\tilde X^{n+1})-\EE(X^n) \geq 0.
\end{cases}
\label{s2m}
\ee
The above relationship, given the two search positions $X^n$ and $\tilde X^{n+1}$, results in $X^{n+1}=\tilde X^{n+1}$ if $\tilde X^{n+1}$ is better than $X^n$, i.e., $\EE(\tilde X^{n+1}) < \EE(X^{n})$, otherwise it interpolates between $X^n$ and $\tilde X^{n+1}$ with a weight proportional to the Boltzmann-Gibbs' measure. Therefore, it differs from the standard simulated annealing update in the case where $\tilde X^{n+1}$ is a worse position. In the latter case, instead of accepting that position with a probability dependent on the Boltzmann-Gibbs' measure, the process moves toward that position by an amount proportional to the same Boltzmann-Gibbs' measure. 

In a continuous setting we can write the update rule in compact form as
\be
x' = x + B_\EE(x\to x+\sigma(t)\xi)\sigma(t)\xi,\qquad B_\EE(x\to x+\sigma(t)\xi)=\min\left\{1,\frac{f^\infty_{\EE}(x+\sigma(t)\xi)}{f^\infty_{\EE}(x)}\right\}.
\label{eq:mc}
\ee  

As a result, the Boltzmann-Gibbs' measure appears in the kinetic equation in the microscopic interaction \eqref{eq:mc} and not as a kernel of the integral describing the time variation.  
The corresponding kinetic equation can be written in weak form as
\be
\frac{\partial}{\partial t}\int_{\RR^d} f(x,t)\phi(x)\,dx = \left\langle \int_{\RR^{d}}(\phi(x')-\phi(x))f(x,t)\,dx\right\rangle,
\label{weaksa3}
\ee
for any regular function $\phi=\phi(x)$. 

Let us now carry on the same asymptotic mean-field approximation leading to the Langevin dynamics \eqref{eq:mfsa} in the case of the kinetic description \eqref{eq:ksa}. The starting point is the same scaling \eqref{eq:scal} and expansions \eqref{eq:expand1}-\eqref{eq:expand2} for small values of $\varepsilon \ll 1$. Now we have
\[
\begin{split}
x'-x &= \sqrt{\varepsilon}\sigma(t)\xi - \left(1-\frac{f^\infty_{\EE}(x+\sqrt{\varepsilon}\sigma(t)\xi)}{f^\infty_{\EE}(x)}\right)\psi(\EE(x+\sqrt{\varepsilon}\sigma(t)\xi)\geq \EE(x))\sqrt{\varepsilon}\sigma(t)\xi\\
&= \sqrt{\varepsilon}\sigma(t)\xi - \varepsilon\sigma(t)\xi\cdot\frac1{T(t)}\nabla_x\EE(x)\psi(\EE(x+\sqrt{\varepsilon}\sigma(t)\xi)\geq \EE(x))\sigma(t)\xi+O(\varepsilon^{3/2}),
\end{split}
\]
and, using the fact that $\xi$ has zero mean, it is immediate to verify that we obtain the same approximated weak form \eqref{eq:exp1} up to $O(\sqrt{\varepsilon})$ as in the classical case.

By the same arguments using the symmetry of $p(\xi)$ and taking $2T(t)=\sigma^2(t)$ we conclude that in the limit $\varepsilon\to 0$ we get the Langevin dynamics \eqref{eq:VlasovSA}.

\begin{remark}^^>
\begin{itemize}
\item Although the two variants of the simulated annealing optimization process here described, at the leading order, lead to the same asymptotic expansion \eqref{eq:exp1} and thus to the Langevin dynamics \eqref{eq:VlasovSA}, the higher order terms in the expansions of the kinetic optimization models are clearly different leading to different behaviours for small but non zero values of the scaling parameters.
\item Further variants of the kinetic simulated annealing dynamics can be considered, for example taking an anisotropic temperature to reduce dimensional dependence in convergence to the global minimum as in CBO^^>\cite{carrillo2019consensus}, or considering the temperature as an independent variable of the system which evolves accordingly to a given dynamic as in parallel tempering^^>\cite{MP92}. These aspects will be the subject of future studies.
\end{itemize}
\end{remark}

\section{Algorithmic aspects and numerical validation}
The kinetic equation \eqref{eq:ksa} can be solved using a direct simulation Monte Carlo method combined with a suitable acceptance-rejection technique which stress the analogies with variable hard-spheres simulations in rarefied gas dynamics^^>\cite{partos13}. 

More precisely, after introducing a time step $\Delta t$ we consider the time discrete form of the scaled kinetic simulated annealing \eqref{eq:ksa} in the form
\be
f^{n+1}(x) = f^n(x)\left(1-\frac{\Delta t}{\varepsilon}\right) + \frac{\Delta t}{\varepsilon} \mathcal{J}^\varepsilon_\EE(f^n)(x),
\label{eq:ksad}
\ee
where $\mathcal{J}^\varepsilon_\EE(f)(x)=\mathcal{L}^\varepsilon_\EE(f)(x)+f(x)$ is the scaled operator accordingly to \eqref{eq:scal}. Now, thanks to \eqref{eq:J}, this operator is nonnegative and can be understood as a probability density function. Therefore, equation \eqref{eq:ksad} for $\Delta t \leq \varepsilon$ defines a convex combination between two probability densities. In the sequel we assume $\Delta t =\varepsilon$, and $\sigma=\sqrt{2T}$ so that the discrete process reduces to
\[
f^{n+1}(x) = \mathcal{J}^\varepsilon_\EE(f^n)(x),
\] 
where sampling from the right hand side corresponds to moving the search point accordingly to the trial position
\[
x'=x+\sqrt{2\varepsilon T}\xi
\]
and accepting or rejecting the position with probability $B_\EE(x\to x')$ defined in \eqref{eq:kernel}.

Here we report in Algorithm \ref{algo:ksa} the time-discrete scheme under the assumption $\Delta t = \varepsilon$ in \eqref{eq:ksad} following^^>\cite{AlPa}. Its Maxwellian counterpart is realized similarly simply replacing the acceptance-rejection routine with the interpolation schedule defined in \eqref{s2m}. The details are given in Algorithm \ref{algo:msa}.

\medskip
\begin{algorithm}[H]
\begin{algorithmic}
\footnotesize
\STATE{Input parameters: $T_0$, $n_t$, $\varepsilon>0$}
\STATE{Initialize trial point $X^{(0)}$} 
\STATE{$T^{(0)}\gets T_0$, $t \gets 0$}
\WHILE{$t < n_t$}
\STATE{Generate $\xi \sim \mc{N}(0,1)$}
\STATE{Compute $\tilde X = X^{(t)}+\sqrt{2\varepsilon T^{(t)}}\xi$} 
\STATE{Evaluate $B^{(t)}_\EE=B_\EE(X^{(t)}\to \tilde X)$ defined in \eqref{eq:kernel}} 
\STATE{Generate $\eta \sim \mc{U}(0,1)$}
\IF{$\eta \leq B^{(t)}_\EE$}
\STATE{$X^{(t+1)}\gets \tilde X$}
\ELSE
\STATE{$X^{(t+1)}\gets X^{(t)}$}
\ENDIF
\STATE{$t\gets t+1$}
\STATE{$T^{(t)}\gets T_0/\log(t+2)$}
\ENDWHILE
\end{algorithmic}
\caption{Kinetic Simulated Annealing (KSA)}\label{algo:ksa}
\end{algorithm}

\begin{algorithm}[H]
\begin{algorithmic}
\footnotesize
\STATE{Input parameters: $T_0$, $n_t$, $\varepsilon>0$}
\STATE{Initialize trial point $X^{(0)}$} 
\STATE{$T^{(0)}\gets T_0$, $t \gets 0$}
\WHILE{$t < n_t$}
\STATE{Generate $\xi \sim \mc{N}(0,1)$}
\STATE{Compute $\tilde X = X^{(t)}+\sqrt{2\varepsilon T^{(t)}}\xi$} 
\STATE{Evaluate $B^{(t)}_\EE=B_\EE(X^{(t)}\to \tilde X)$ defined in \eqref{eq:kernel}} 
\STATE{$X^{(t+1)}\gets X^{(t)}+B^{(t)}_\EE(\tilde X-X^{(t)})$}
\STATE{$t\gets t+1$}
\STATE{$T^{(t)}\gets T_0/\log(t+2)$}
\ENDWHILE
\end{algorithmic}
\caption{Maxwellian Simulated Annealing (MSA)}\label{algo:msa}
\end{algorithm}

\medskip

We would like to emphasize that our aim is not to show the performance of the methods on a multitude of test functions, as our results are mainly of theoretical interest, but to compare kinetic and simulated Maxwellian annealing and to validate the corresponding mean-field Langevin limit. The latter result may be particularly relevant in applications, as it shows that gradient-based results can be obtained with gradient-free methods under appropriate parameter scaling. To this end, we restrict our results to the Ackley function, a prototype test function in global optimization that has previously been adopted as a reference for mean-field optimization dynamics^^>\cite{pinnau2017consensus,carrillo2019consensus}. 

\begin{figure}[htb]
\includegraphics[scale=0.45]{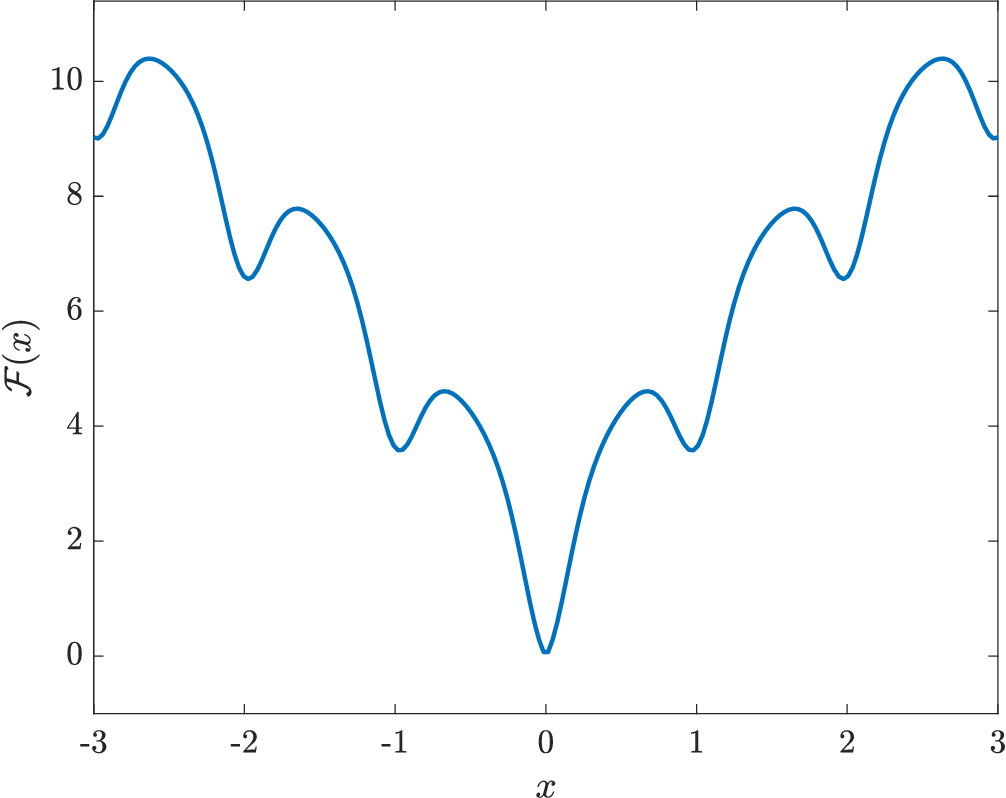}
\includegraphics[scale=0.45]{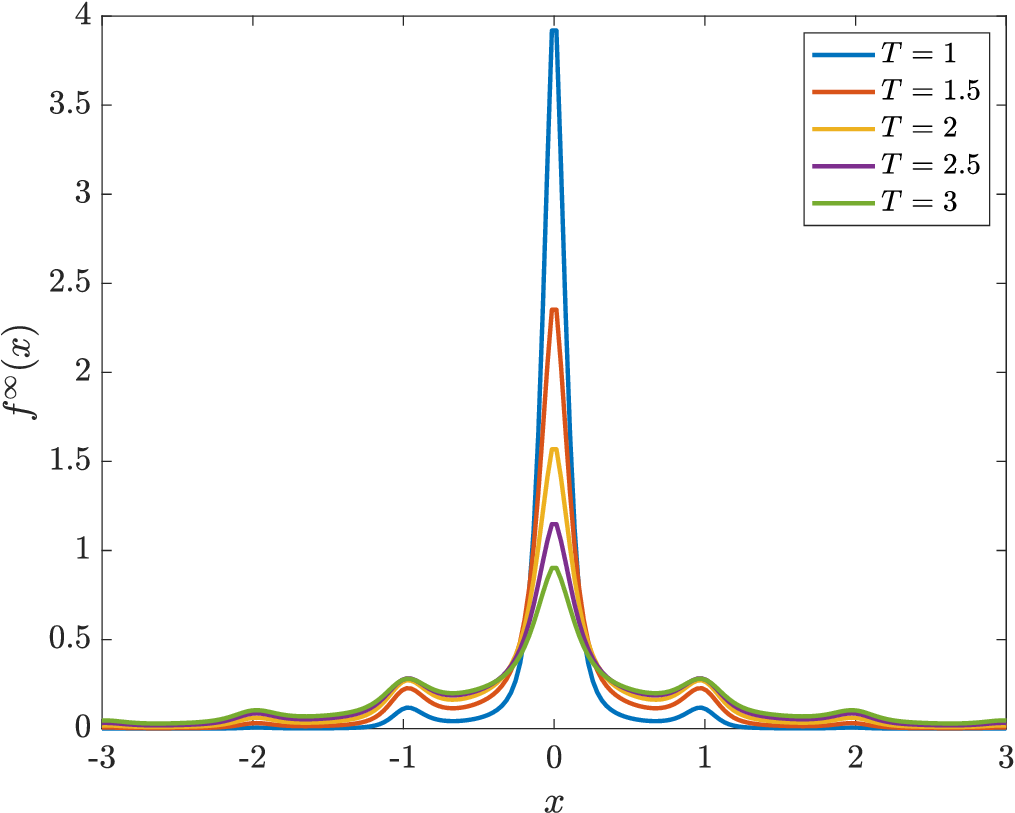}
\caption{The prototype Ackley function (left) and the corresponding steady states (right) given by the Gibbs measure \eqref{eq:Gibbs} for various values of the control temperature.}
\label{fg:fig0}
\end{figure}

\begin{figure}[htb]
\includegraphics[scale=0.45]{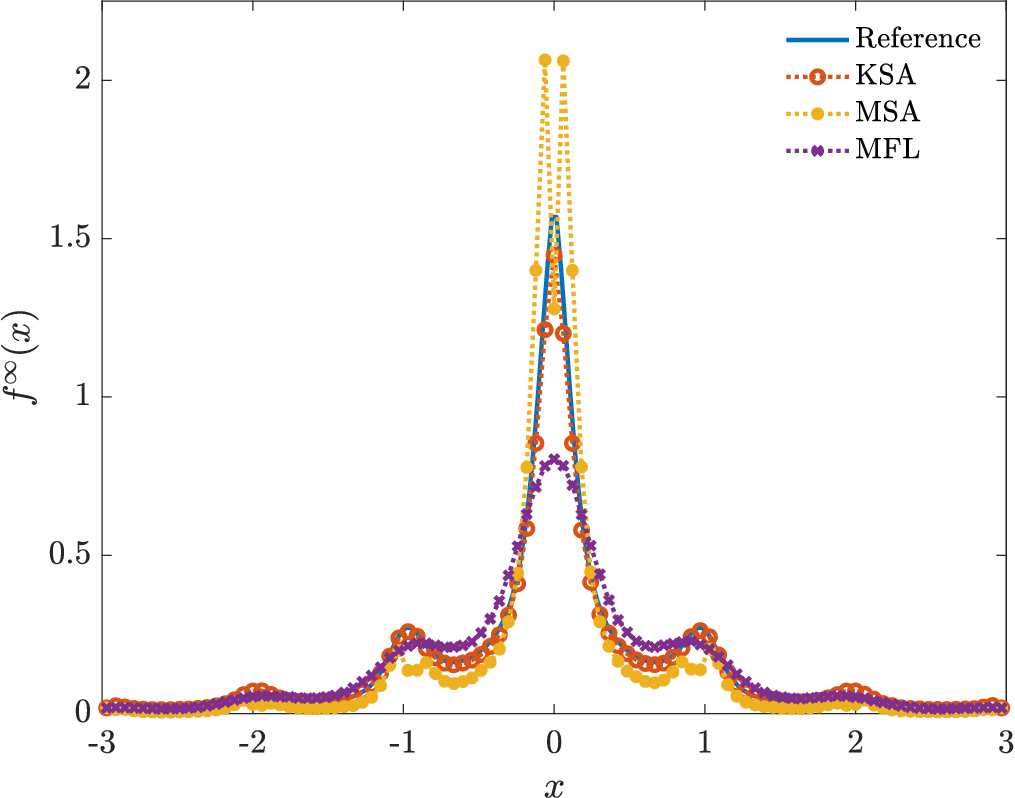}
\includegraphics[scale=0.45]{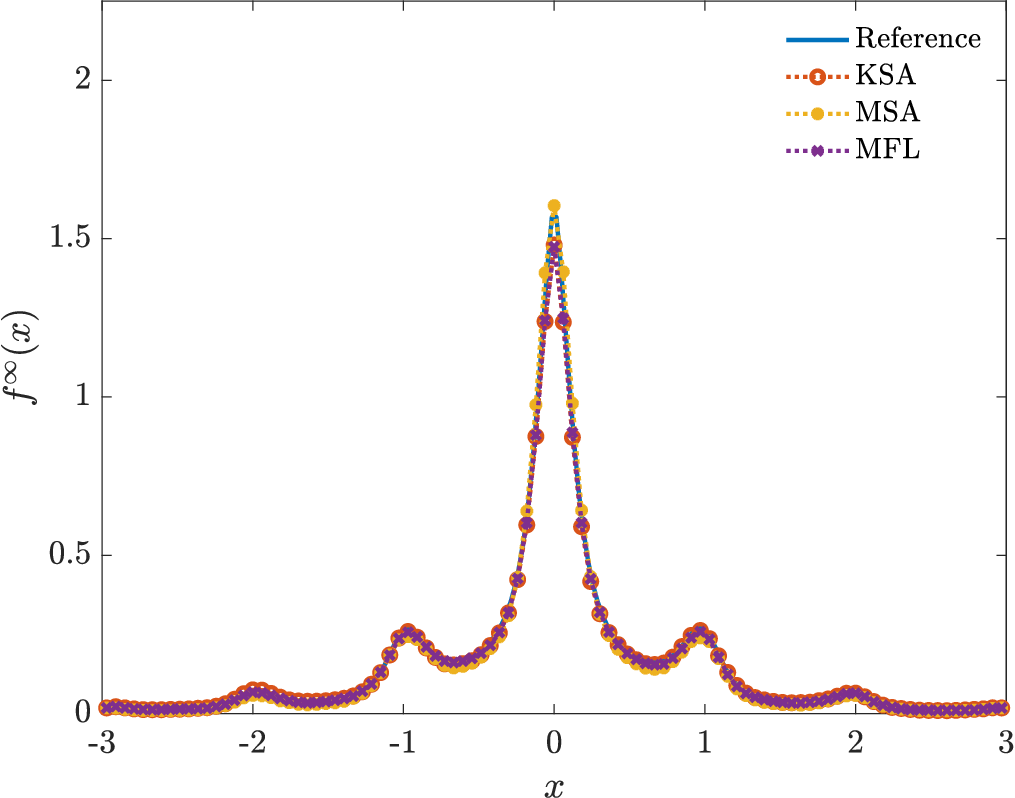}
\includegraphics[scale=0.45]{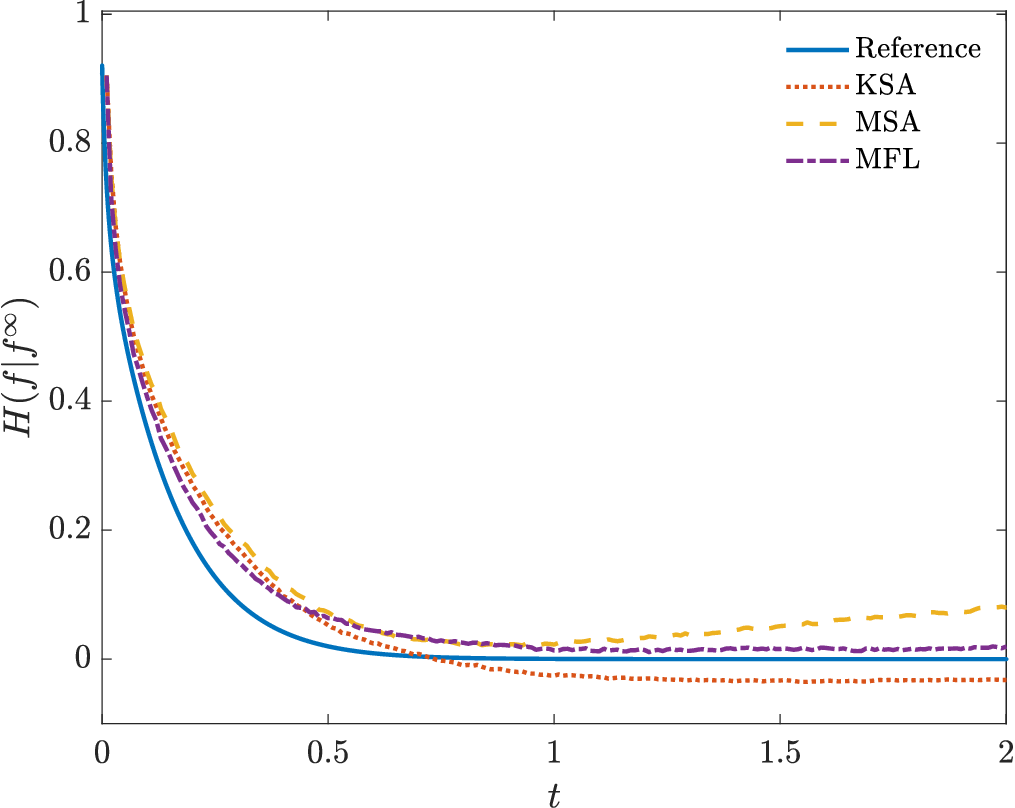}
\includegraphics[scale=0.45]{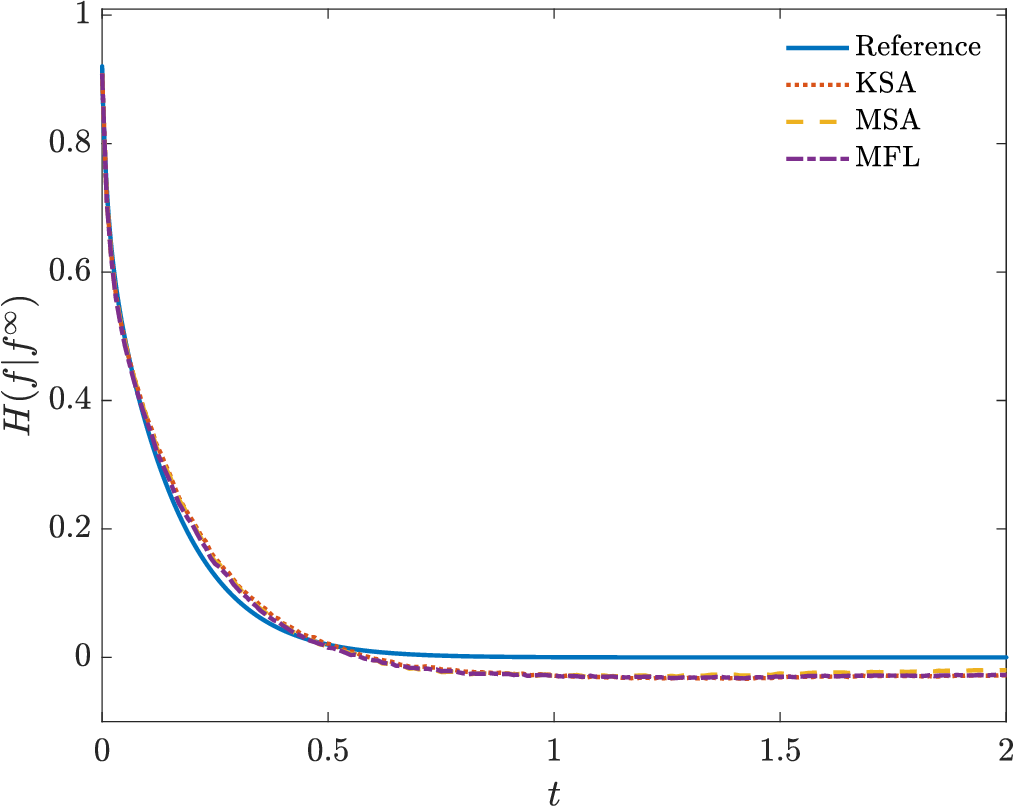}
\caption{Solution of KSA and MSA for a fixed control temperature $T=2$ for $\varepsilon=0.01$ (left) and $\varepsilon=0.0001$ (right). On the top the probability density at final time $t=2$, on the bottom relative entropies along the simulation. As a reference we also report the mean-field (Reference) and the stochastic Langevin dynamic (MFL) results. All plots have been obtained averaging over $N=5\times 10^4$ runs.}
\label{fg:fig1}
\end{figure}

\begin{figure}[H]
\includegraphics[scale=0.45]{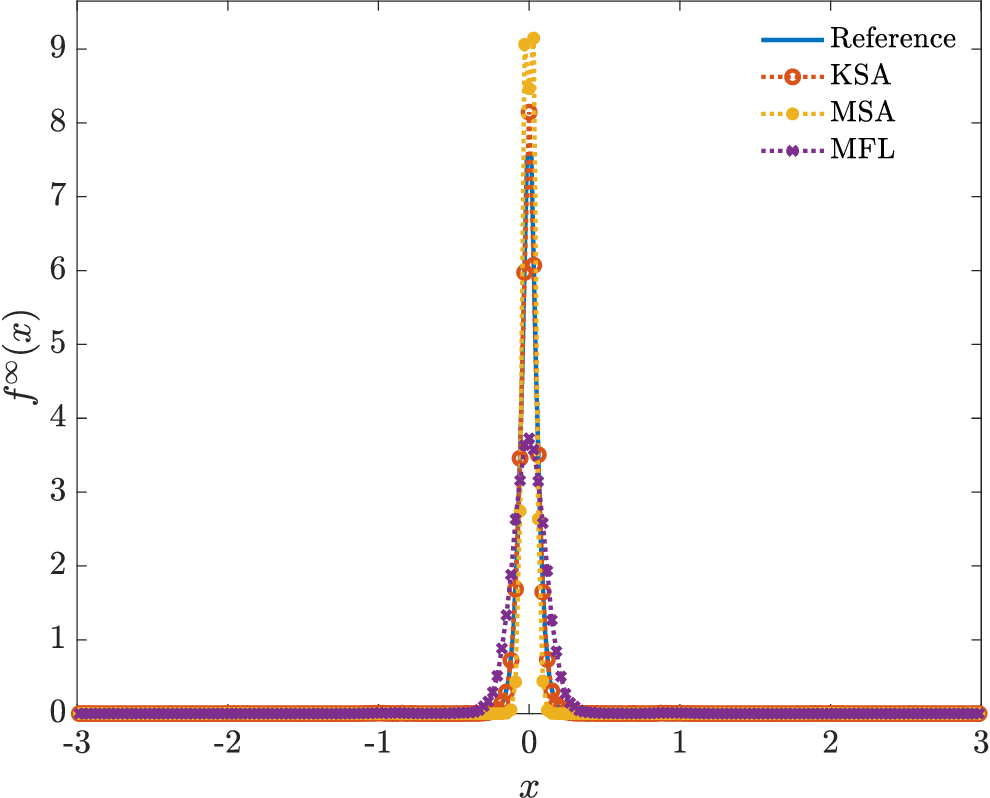}
\includegraphics[scale=0.45]{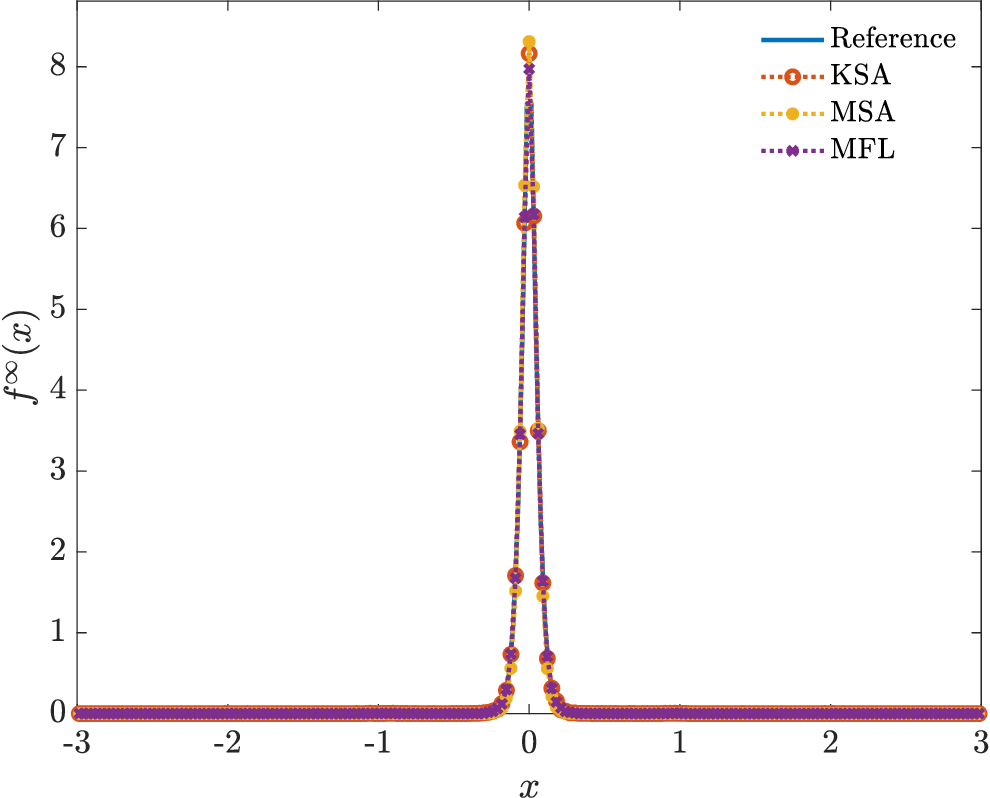}
\includegraphics[scale=0.45]{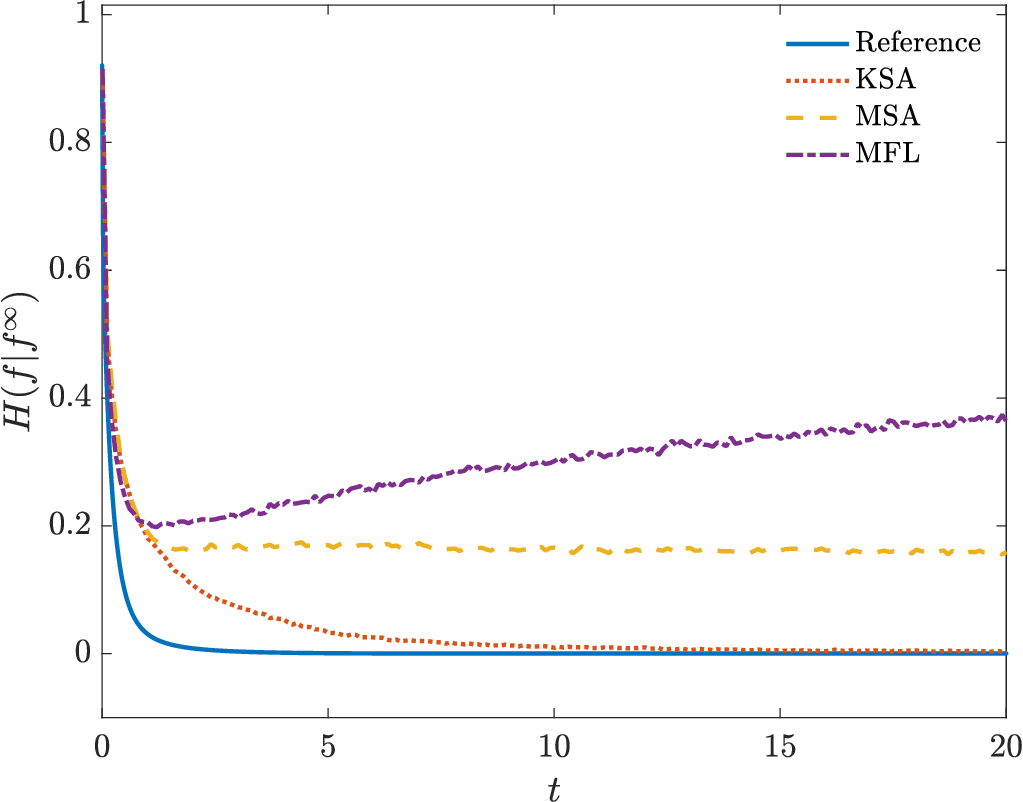}
\includegraphics[scale=0.45]{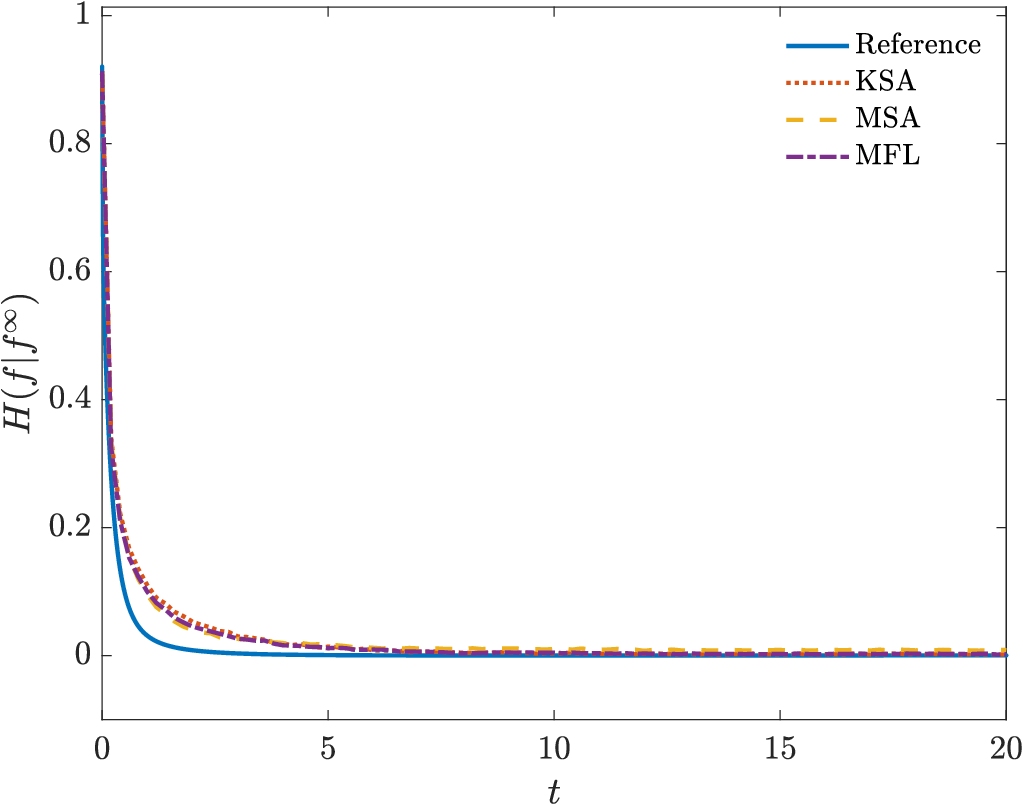}
\caption{Solution of KSA and MSA for a time-dependent control temperature $T(t)=2\log(2)/\log(2+t)$ for $\varepsilon=0.01$ (left) and $\varepsilon=0.0001$ (right). On the top the probability density at final time $t=20$, on the bottom relative entropies along the simulation. As a reference we also report the mean-field (Reference) and the stochastic Langevin dynamic (MFL) results. All plots have been obtained averaging over $N=50.000$ runs.}
\label{fg:fig2}
\end{figure}

In addition to the KSA and MSA methods, we considered as comparisons the solutions obtained by solving the Langevin mean-field dynamics \eqref{eq:mfsa} with a structure-preserving finite difference method that is able to correctly capture the long-term behaviour of the system^^>\cite{PZ18}, and by solving the corresponding SDE \eqref{eq:csa} with the Euler-Maruyama method. The latter approach, referred to as MFL and similar in practice to SGD, was applied by taking $\Delta t=\varepsilon$. All test cases have been performed in a simple one-dimensional setting starting from a uniform distribution in $[-3,3]$ as initial data. The results are reported in Figure \ref{fg:fig1} for a fixed temperature $T=2$ and in Figure \ref{fg:fig2} for a time dependent temperature $T(t)=2\log(2)/\log(t+2)$.

It can be observed that in the constant-temperature case of Figure \ref{fg:fig1} the KSA method converges towards the Gibbs steady state and preserves the decay of relative entropy independently of the scaling term $\varepsilon$. 
In contrast, the MSA method, as well as the MFL method, have different behaviours for $\varepsilon=0.01$. In particular, for the MSA method, the relative entropy does not decrease and the steady state depends on the rescaling parameter, which presumably makes a rigorous mathematical analysis of the long time behavior of the method rather difficult.  As expected, all methods for smaller values of the rescaling parameter, $\varepsilon=0.0001$, converge to the same mean-field solution of the Langevin dynamics. 

Next, we considered a time-dependent control $T(t)=2\log(2)/\log(t+2)$ which should guarantee for KSA convergence to the global minimum. In Figure \ref{fg:fig2} we report the results. The decay of the relative entropy for KSA independently from the scaling is confirmed by the numerical experiments. Furthermore, we observe how MSA, as well as MFL, recover the entropy dissipation for small values of $\varepsilon$. In particular, all methods converges to the global minimum, with MSA having a lower variance around the minimum and MFL a larger one. Again, as expected, all methods converge to the mean-field solution of the Langevin dynamic when $\varepsilon$ tends to zero.    

\section{Concluding remarks}
We have introduced novel formulations of the simulated annealing method for global optimization based on linear kinetic equations. This allows to highlight analogies between classical entropy-based approaches in studying the trend towards equilibrium of kinetic equations and the convergence of the method to the global minimum of the objective function.

By adopting a continuous setting based on kinetic partial differential equations, we explored the relationships with Fokker-Planck equations that describe the mean-field Langevin dynamic. Specifically, we demonstrated how to formally derive the corresponding mean-field model, which utilizes gradient information, taking a suitable scaling limit of the gradient-free linear kinetic model that describes simulated annealing. This derivation makes it possible to construct a theoretical link between gradient-based and gradient-free techniques by exploiting the analogies with the grazing collision limit of the Boltzmann equation.  

We extended this methodology to other types of kinetic equations that describe variations of the simulated annealing process, avoiding the acceptance-rejection process. Numerical evidence supporting such asymptotic behavior was also illustrated through some simple simulation examples. The purpose of these examples is to validate the theoretical results presented here, leaving the possible construction and implementation of new high-dimensional efficient optimization methods based on these ideas to further research.

From a mathematical perspective, it is important to note that several challenging questions remain open. For instance, estimating the rate of convergence to equilibrium or analyzing the convergence properties of the Maxwellian variant introduced here.

\section*{Acknowledgment}
This work has been written within the activities of GNCS group of INdAM (Italian National Institute of High Mathematics). The research has been supported by the Royal Society under the Wolfson Fellowship ``Uncertainty quantification, data-driven simulations and learning of multiscale complex systems governed by PDEs". The partial support by ICSC -- Centro Nazionale di Ricerca in High Performance Computing, Big Data and Quantum Computing, funded by European Union -- NextGenerationEU and by MIUR-PRIN Project 2022, No. 2022KKJP4X ``Advanced numerical methods for time dependent parametric partial differential equations with applications" is also acknowledged.

\bibliographystyle{abbrv}
\bibliography{biblio}
\end{document}